\documentclass{amsart}
\usepackage{amssymb}
\usepackage{mathrsfs}
\usepackage{bbm}
\usepackage{cite}
\usepackage{ifthen}

\newtheorem{theorem}{Theorem}[section]
\newtheorem{lemma}[theorem]{Lemma}
\newtheorem{cor}[theorem]{Corollary}
\newtheorem{prop}[theorem]{Proposition}

\theoremstyle{remark}
\newtheorem{rem}[theorem]{Remark}

\theoremstyle{definition}
\newtheorem{mdef}[theorem]{Definition}

\newcommand*{\R}{\ensuremath{\mathbb{R}}}
\newcommand*{\N}{\ensuremath{\mathbb{N}}}
\newcommand*{\T}{\ensuremath{\mathbb{T}}}
\newcommand*{\ls}{\lesssim}
\newcommand*{\gs}{\gtrsim}
\newcommand*{\jbr}[1]{\langle #1 \rangle}
\newcommand*{\eps}{\varepsilon}
\newcommand*{\del}{\delta}
\newcommand*{\la}{\lambda}
\newcommand*{\al}{\alpha}
\newcommand*{\lf}{\sigma}
\newcommand*{\ftr}[1]{\widehat{#1}}
\newcommand*{\alftr}{\ensuremath{\mathcal{F}}}
\newcommand*{\xsb}[4][0]{
  \ifthenelse{\equal{#3}{0}}{
  \ensuremath{X^{#2,#4}_{#1}}
  }{
  \ensuremath{X^{#2,#3,#4}_{#1}}  
}}
\newcommand*{\embeds}{\hookrightarrow}
\newcommand*{\ntf}{\ensuremath{\tilde{\mathcal{S}}}}
\newcommand*{\tf}{\ensuremath{\mathcal{S}}}
\newcommand*{\z}{\mu}

\DeclareMathOperator{\sign}{sign}
\DeclareMathOperator{\med}{med}

\title[Well-posedness for KP II and generalisations]
{Well-posedness for the Kadomtsev-Petviashvili II equation
and generalisations}
\author{Martin Hadac}

\begin{document}

\bibliographystyle{hplain}

\begin{abstract}
  We show the local in time well-posedness of the Cauchy problem
  for the Kadomtsev-Petviashvili II equation for initial data
  in the non-isotropic Sobolev space $H^{s_1,s_2}(\R^2)$
  with $s_1>-\frac12$ and $s_2\geq 0$.
  On the $H^{s_1,0}(\R^2)$ scale this result includes the full subcritical
  range without any additional low frequency assumption on the
  initial data.
  More generally, we prove the local in time well-posedness 
  of the Cauchy problem for the
  following generalisation of the KP II equation
  \(
    (u_t - |D_x|^\al u_x + (u^2)_x)_x + u_{yy} = 0,
    \quad
    u(0) = u_0
  \)
  for $\frac43<\al\leq 6$, $s_1>\max(1-\frac34 \al,\frac14-\frac38 \al)$, 
  $s_2\geq 0$ and $u_0\in H^{s_1,s_2}(\R^2)$.
  We deduce global well-posedness for $s_1\geq 0$, $s_2=0$
  and real valued initial data.
\end{abstract}

\maketitle

\section{Introduction}
\
In this paper we study the Cauchy problem for the
Kadomtsev-Petviashvili II equation
\begin{equation}
  \label{eq:DGKP_kpeq}
  (u_t + u_{xxx} + (u^2)_x)_x + u_{yy} = 0
  \quad \text{in $\R^3$},\quad
  u(0) = u_0.
\end{equation}
The Kadomtsev-Petviashvili II equation (as well as the 
Kadomtsev-Petviashvili I equation 
$(u_t + u_{xxx} + (u^2)_x)_x - u_{yy} = 0$)
are two-dimensional extensions of the Korteweg-de-Vries equation,
see \cite{kadomtsev:70}.
More generally, we will consider the following dispersion generalised
Kadomtsev-Petviashvili II type equation
\begin{equation}
  \label{eq:DGKP_dgkpeq}
  (u_t - |D_x|^\al u_x + (u^2)_x)_x + u_{yy} = 0
  \quad \text{in $\R^3$},\quad
  u(0) = u_0
\end{equation}
where $\frac43 < \al \leq 6$.
Here $|D_x|^\al$ is the Fourier multiplier operator with multiplier
$|\xi|^\al$.
We consider initial values $u_0$ in the non-isotropic Sobolev spaces 
\begin{equation}
  \label{eq:DGKP_defsob}
  H^{s_1,s_2}(\R^2) := \{ u_0\in \tf'(\R^2) \mid \|u_0\|_{H^{s_1,s_2}} :=
  \| \jbr{\xi}^{s_1} \jbr{\eta}^{s_2} \ftr{u_0} \|_{L^2_{\xi\eta}} < \infty \}.
\end{equation}
Note that the case $\al=2$ of \eqref{eq:DGKP_dgkpeq} is 
\eqref{eq:DGKP_kpeq} whereas the case $\al=4$ is known as 
fifth order KP II equation
\begin{equation}
  \label{eq:DGKP_fokpeq}
  (u_t - u_{xxxxx} + (u^2)_x)_x + u_{yy} = 0
  \quad \text{in $\R^3$},\quad
  u(0) = u_0.  
\end{equation}

We are interested in low regularity well-posedness of 
\eqref{eq:DGKP_dgkpeq}.
By using refined Fourier restriction norm spaces we will prove
new bilinear estimates which allow us to apply the contraction
mapping principle.

In the seminal work \cite{bourgain:93} Bourgain shows the 
(global) well-posedness of \eqref{eq:DGKP_kpeq}
(on $\T^2$ rather than on $\R^2$) with initial values in $L^2$,
i. e. for $s_1=s_2=0$.
This result has been improved afterwards by 
Takaoka and Tzvetkov \cite{takaoka:01} and 
Isaza and Mej\'ia \cite{isaza:01} to the local in time well-posedness
of \eqref{eq:DGKP_kpeq} for $s_1>-\frac13$ and $s_2\geq 0$.
(For previous results see also \cite{tzvetkov:99},
\cite{tzvetkov:00}, \cite{takaoka:01_2}.)
In \cite{takaoka:00} Takaoka shows local well-posedness
for $s_1>-\frac12$, $s_2=0$, but only if the additional low frequency
condition $|D_x|^{-\frac12+\eps} u_0\in L^2$
(with suitably chosen $\eps$) is imposed on the initial value.

For the fifth order KP-II equation \eqref{eq:DGKP_fokpeq} local 
well-posedness was shown by Saut and Tzvetkov (see 
\cite{saut:99,saut:00}) for $s_1\geq -\frac14$ and $s_2\geq 0$.
(Note that the equation considered in \cite{saut:99,saut:00}
is slightly more general than \eqref{eq:DGKP_fokpeq} because
it also contains the third order term.)
Very recently, Isaza, L\'opez and Mej\'ia \cite{isaza:06}
improved the local well-posedness result to $s_1>-\frac54$
and $s_2\geq 0$. (These authors also show global well-posedness
of \eqref{eq:DGKP_fokpeq} for $s_1>-\frac47$ and $s_2=0$.) 

For general $\al \in (\frac43,6]$ I\'orio and Nunes \cite{iorio:98}
showed the local well-posedness for initial values $u_0$ in the
isotropic Sobolev space $H^s(\R^2)$, $s>2$, with the additional
low frequency condition $\partial_x^{-1} u_0\in H^s(\R^2)$
using parabolic regularization. 
Let us note that they consider much more general equations
and do not use the dispersive structure of the equation.

For a recent result concerning the so called mass constraint
property for solutions of equations of type \eqref{eq:DGKP_dgkpeq}
see \cite{saut:06}.

Our main result for equation \eqref{eq:DGKP_kpeq} is the following
\begin{theorem}
  \label{th:DGKP_mainkp}
  Let $s_1>-\frac12$ and $s_2\geq 0$.
  For $R>0$ there exists $T=T(R)>0$ and a Banach space
  $X_T \embeds C([-T,T];H^{s_1,s_2}(\R^2))$, such that for every 
  $u_0\in B_R 
  := \{u_0\in H^{s_1,s_2}(\R^2) \mid \|u_0\|_{H^{s_1,s_2}(\R^2)} < R \}$ 
  there is exactly one solution $u$ of equation 
  \eqref{eq:DGKP_kpeq} in $X_T$.
  Furthermore the mapping 
  $F_R:B_R\to X_T, u_0\mapsto u$ is analytic.
\end{theorem}
More generally, we will show the following theorem concerning
equation \eqref{eq:DGKP_dgkpeq}
\begin{theorem}
  \label{th:DGKP_maindgkp}
  Let $\frac43<\al\leq 6$, $s_1>\max(1-\frac34 \al,\frac14-\frac38 \al)$
  and $s_2\geq 0$.
  For $R>0$ there exists $T=T(R)>0$ and a Banach space
  $X_T \embeds C([-T,T];H^{s_1,s_2}(\R^2))$, such that for every 
  $u_0\in B_R 
  := \{u_0\in H^{s_1,s_2}(\R^2) \mid \|u_0\|_{H^{s_1,s_2}(\R^2)} < R \}$ 
  there is exactly one solution $u$ of equation 
  \eqref{eq:DGKP_dgkpeq} in $X_T$.
  Furthermore the mapping 
  $F_R:B_R\to X_T, u_0\mapsto u$ is analytic.
\end{theorem}

\begin{rem}
  If we (formally) apply the operator $\partial_x^{-1}$
  to equation \eqref{eq:DGKP_dgkpeq} and use Duhamel's formula,
  equation \eqref{eq:DGKP_dgkpeq} is (for suitable $u$) equivalent
  to the integral equation
  \begin{equation}
    \label{eq:DGKP_dgkpint}
    u(t) = U_\al(t)u_0 - \int_0^t U_\al(t-t')\partial_x(u(t')^2) dt'
  \end{equation}
  where $U_\al$ is the unitary group on $H^{s_1,s_2}(\R^2)$ defined by
  \begin{equation}
    \label{eq:DGKP_linsol}
    \ftr{U_\al(t)u_0}(\xi,\eta) 
    := \exp(it(\xi|\xi|^\al - \frac{\eta^2}{\xi}))\ftr{u_0}(\xi,\eta).
  \end{equation}
  We define a solution of \eqref{eq:DGKP_dgkpeq} 
  in $X_T$ (for $T\leq 1$) to be a solution of the operator equation
  \begin{equation}
    \label{eq:DGKP_opeq}
    u(t) = \psi(t) U_\al(t) u_0 - \Gamma_T (u,u)(t), \quad t\in [-T,T] 
  \end{equation}
  where $\Gamma_T$ is the bilinear operator on $X_T$ defined
  for smooth $u_1$, $u_2$ by
  \begin{equation}
    \label{eq:DGKP_defgam}
    \Gamma_T(u_1,u_2)(t) := 
    \psi_T(t) \int_0^t U_\al(t-t') \partial_x(u_1 u_2)(t') dt'
  \end{equation}
  and $\psi\in C^\infty_0(\R)$ is a cut-off function with 
  $\psi(t)=1$ for $|t|\leq 1$
  and $\psi(t)=0$ for $|t|\geq 2$.
  Furthermore $\psi_T(t) = \psi(t/T)$.
\end{rem}

\begin{rem}
  In the particular case $\al=4$ of the fifth order KP II equation 
  Theorem~\ref{th:DGKP_maindgkp} shows the local 
  well-posedness of \eqref{eq:DGKP_fokpeq} for 
  $s_1>-\frac54$ and $s_2\geq 0$.
  We therefore get a local well-posedness result for the same class
  of initial data as Isaza, L\'opez and Mej\'ia in \cite{isaza:06}.
  Note though that the spaces $X_T$ where the local well-posedness
  result of Theorem~\ref{th:DGKP_maindgkp} holds true are different
  from those used in \cite{isaza:06} (see Remark~\ref{rm:DGKP_xtsp}).
\end{rem}

\begin{rem}
  Let us note that if $u$ is a solution of \eqref{eq:DGKP_dgkpeq}
  then so is 
  \begin{equation*}
    u_\la(t,x,y) = \la^\al u(\la^{\al+1} t, \la x, \la^{\frac\al2 + 1} y). 
  \end{equation*}
  Considering the homogeneous Sobolev norm
  \begin{equation*}
    \|u_0\|_{\dot{H}^{s_1,s_2}} :=
    \| |\xi|^{s_1} |\eta|^{s_2} \ftr{u_0} \|_{L^2_{\xi,\eta}}      
  \end{equation*}
  we get 
  $\|u_\la(0,\cdot,\cdot)\|_{\dot{H}^{s_1,s_2}}
  =\la^{\frac34\al-1+s_1+(1+\frac\al2)s_2}
  \|u(0,\cdot,\cdot)\|_{\dot{H}^{s_1,s_2}}$.
  This argument suggests that we get ill-posedness for 
  $s_1+(1+\frac\al2)s_2 < 1-\frac34\al$.
  Note that for $\frac43<\al\leq 2$ and $s_2=0$ we
  reach the critical value $1-\frac34\al$ of $s_1$,
  except for the endpoint.
  For $\al>2$ though we have that 
  $\frac14-\frac38\al > 1-\frac34 \al$, so that we do not
  reach the scaling limit in this case.
\end{rem}

By combining the local well-posedness result of 
Theorem~\ref{th:DGKP_maindgkp} with the conservation
of the $L^2$-norm which holds for real valued
solutions of \eqref{eq:DGKP_opeq} we get the following 
global result, where $H^{s_1,0}(\R^2;\R)$ denotes
the subspace of all real valued functions in $H^{s_1,0}(\R^2)$.

\begin{theorem}
  \label{th:DGKP_mainglob}
  Let $\frac43<\al\leq 6$, $s_1\geq 0$ and $T>0$.
  Then there exists a Banach space 
  $X_T \embeds C([-T,T];H^{s_1,0}(\R^2;\R))$, such that for every 
  $u_0\in H^{s_1,0}(\R^2;\R)$ there is exactly one solution $u$ of 
  equation \eqref{eq:DGKP_dgkpeq} in $X_T$.
\end{theorem}

Let us fix some notation we use throughout the paper:
\begin{itemize}
\item For $\xi\in\R$ let $\jbr{\xi} := (1+|\xi|^2)^{\frac12}$.
\item For $u\in \tf'(\R^n)$ the Fourier transformation
  of $u$ in $\R^n$ is denoted by $\ftr{u}$ or $\alftr{u}$.
  A partial Fourier transformation with respect to some of 
  the $n$ variables, is denoted for example by 
  $\alftr_1$ for the Fourier transformation in the first variable, etc.
\item $\z=(\tau,\xi,\eta)\in \R^3$ always denotes the Fourier variable
  dual to $(t,x,y)$.
\item For $\z=(\tau,\xi,\eta)$ let
  $\la := \la(\z) := \tau - \xi|\xi|^\al + \frac{\eta^2}{\xi}$.
  If there are two frequency variables $\z$ and $\z_1$ we write
  for short $\la_1 := \la(\z_1)$, $\la_2 := \la(\z-\z_1)$ and
  $|\la_{\max}|=\max(|\la|,|\la_1|,|\la_2|)$.
  Let also $|\xi_{\max}| := \max (|\xi|, |\xi_1|, |\xi-\xi_1|)$.
  Let $|\xi_{\min}|$ and $|\xi_{\med}|$ be defined analogously.
\item $A\ls B$ means that there is a (harmless) constant
  $C$ such that $A\leq CB$.
\item For $X$ and $Y$ Banach spaces $X\embeds Y$ means that there
  is a continuous embedding from $X$ into $Y$.
  Furthermore $C_b(\R;X)$ denotes the space of all continuous and
  bounded functions $f:\R\to X$ with the $\sup$-norm.
\end{itemize}

The author would like to thank S. Herr and H. Koch for valuable
discussions and suggestions on the subject.

\section{Definition of the solution spaces}
\label{se:DGKP_defxsb}
\begin{mdef}
  Let us consider the following space of test functions
  \begin{equation}
    \label{eq:DGKP_tfdef}
    \ntf := \{ \phi \in \tf(\R^3) \mid \partial_\xi^k \ftr{\phi}(\tau,0,\eta)=0
    \ \forall k\in \N_0\ \forall (\tau,\eta)\in \R^2 \}.
  \end{equation}
  For $s_1,s_2,b\in \R$ , $\lf\geq 0$ and $\phi\in\ntf$ define
  \begin{equation}
    \label{eq:DGKP_xsbndef}
    \|\phi\|_{\xsb[\lf]{s_1}{s_2}b} 
    := \| |\xi|^{-\lf} \jbr{\xi}^{s_1+\lf} \jbr{\eta}^{s_2}
    \jbr{\la}^b \ftr{\phi} \|_{L^2_\z}.
  \end{equation}
  Let $\xsb[\lf]{s_1}{s_2}b$ be the completion
  of $\ntf$ with respect to the norm \eqref{eq:DGKP_xsbndef}. 
\end{mdef}

Functions in $\ntf$ have the property that for every $k\in \N_0$ there is
a $C_k>0$ such that $|\ftr{\phi}(\tau,\xi,\eta)|\leq C_k |\xi|^k$
for all $(\tau,\xi,\eta)\in \R^3$.
This property ensures that the right hand side of~\eqref{eq:DGKP_xsbndef} 
is well defined in spite of the singularity along $\xi=0$ in the term 
$\la=\tau-\xi|\xi|^\al+\frac{\eta^2}\xi$ 
and the factor $|\xi|^{-\lf}$.

\begin{rem}
  At least for $b>-\frac12-\sigma$,
  we can identify $\xsb[\lf]{s_1}{s_2}b$ with the subspace of
  tempered distributions $u$ on $\R^3$ 
  such that $\ftr{u}$ is a regular distribution and
  $|\xi|^{-\lf} \jbr{\xi}^{s_1+\lf} \jbr{\eta}^{s_2}\jbr{\la}^b \ftr{u}\in L^2$. 
\end{rem}

\begin{rem}
  If $s_2=0$ we write for short $\xsb[\lf]{s_1}0b$ instead of
  $X^{s_1,s_2,b}_{\lf}$.
\end{rem}

\begin{rem}
  The spaces $\xsb[\lf]{s_1}{s_2}b$ are modifications of 
  spaces first used by Bourgain \cite{bourgain:93_2,bourgain:93_3}
  in the context of the KdV and Schr\"odinger equations.
\end{rem}

We have the following well-known linear estimates
\begin{prop}
  \label{pr:DGKP_linesth}
  For $b\geq 0$ and $s_1,s_2\in \R$ we have
  \begin{equation}
    \label{eq:DGKP_linestiv}
    \| \psi U_\al (t) u_0 \|_{\xsb{s_1}{s_2}b} \ls \|u_0\|_{H^{s_1,s_2}(\R^2)}
  \end{equation}
\end{prop}
\begin{proof}
  See for example \cite{ginibre:96}.
\end{proof}

\begin{prop}
  \label{pr:DGKP_linesti}
  For $-\frac12 < b' \leq 0\leq b\leq b'+1$, $T\leq 1$ and 
  $s_1,s_2\in \R$ we have
  \begin{equation}
    \label{eq:DGKP_linestrc}
    \| \psi_T \int_0^t U_\al(t - t') F(t') dt' \|_{\xsb[\lf]{s_1}{s_2}b}
    \ls T^{1-(b-b')} \|F\|_{\xsb[\lf]{s_1}{s_2}{b'}}.
  \end{equation}
\end{prop}
\begin{proof}
  For $\lf=0$ see \cite{ginibre:96}.
  For $\lf\neq 0$ consider the operator $I_\lf$ defined for $u\in \ntf$ by 
  \(
  (\alftr_2 I_\lf u)(t,\xi,y) 
  = (\frac{\jbr{\xi}}{|\xi|})^\lf \alftr_2 u(t,\xi,y).
  \)
  Then $I_\lf : \xsb[\lf]{s_1}{s_2}b \to \xsb{s_1}{s_2}b$ is an
  isometric isomorphism.
  Therefore we have
  \begin{align*}
    \| \psi_T \int_0^{t} U_\al(t - t') F(t') dt'\|_{\xsb[\lf]{s_1}{s_2}{b'}} 
    & = \| I_\lf \psi_T\int_0^{t}U_\al(t-t')F(t') dt'\|_{\xsb{s_1}{s_2}{b'}}\\
    & = \| \psi_T\int_0^{t}U_\al(t-t')I_\lf F(t') dt'\|_{\xsb{s_1}{s_2}{b'}}\\
    & \ls T^{1-(b-b')} \|I_\lf F\|_{\xsb{s_1}{s_2}b}\\
    & = T^{1-(b-b')} \|F\|_{\xsb[\lf]{s_1}{s_2}b}.
  \end{align*}
\end{proof}

We have the following well-known embedding result for the 
$\xsb[\lf]{s_1}{s_2}{b}$-spaces
\begin{prop}
  Let $s_1,s_2\in \R$, $\lf\geq 0$ and $b>\frac12$. Then
  \begin{equation*}
    \xsb[\lf]{s_1}{s_2}b \embeds C_b(\R;H^{s_1,s_2}(\R^2)).
  \end{equation*}
\end{prop}
\begin{proof}
  For $\lf=0$ see \cite{ginibre:96}.
  But for $\lf>0$ we have that
  \(
  \|u\|_{\xsb{s_1}{s_2}b} \leq \|u\|_{\xsb[\lf]{s_1}{s_2}b}.
  \)
\end{proof}

\begin{mdef}
  Let $X\embeds C_b(\R;H^{s_1,s_2}(\R^2))$.
  Then we define the \emph{restriction norm space}
  $X_T := \{ u\big|_{[-T,T]} \mid u\in X \}$ with norm
  \begin{equation*}
    \| u \|_{X_T} = \inf \{ \|\tilde{u}\|_X \mid \tilde{u}\big|_{[-T,T]} = u \}.
  \end{equation*}
  Then $X_T\embeds C([-T,T];H^{s_1,s_2}(\R^2))$.
\end{mdef}

\section{Strichartz and refined Strichartz estimates}
\
Exactly as in the case $\al=2$ (see \cite{saut:93}) we show the
following Strichartz' estimates for the solution of the linear
equation.
For the convenience of the reader we will give the full proof here.
\begin{theorem}
  Let $2<q\leq \infty$, $\frac1r+\frac1q=\frac12$ 
  and $\gamma := (1-\frac2r)(\frac12-\frac\al 4)$.
  Then we have
  \begin{equation}
    \label{eq:DGKP_linStrich}
    \||D_x|^{-\gamma} U_\al (t)u_0\|_{L^q_t L^r_{xy}} \ls \|u_0\|_{L^2_{xy}}.
  \end{equation}
\end{theorem}
\begin{proof}
  Let $\theta := 1-\frac\al 2$ and $\kappa\in \R$.
  As $\theta<1$, we have that
  \[ 
  m_t(\xi,\eta) 
  := |\xi|^{-\theta+i\kappa}e^{it(\xi|\xi|^\al -\frac{\eta^2}\xi)}\in \tf'(\R^2)
  \] 
  for every $t\in \R$.
  Therefore we have for $u_0\in \tf(\R^2)$ and $t\in \R$
  \begin{equation*}
    |D_x|^{-\theta+i\kappa} U_\al(t) u_0 
    = \alftr^{-1}(m_t \ftr{u_0}) 
    =  \alftr^{-1}(m_t) * u_0 \in \tf'(\R^2).
  \end{equation*}
  For $\del_1, \del_2 > 0$ let us define
  \(
    m^{\del_1,\del_2}_t(\xi,\eta) := e^{-\del_1 \xi^2 -\del_2 \eta^2} m_t(\xi,\eta).
  \)
  Then by the theorem of dominated convergence we have that
  $\lim_{\del_1,\del_2\to 0+} m^{\del_1,\del_2}_t = m_t$ in $\tf'(\R^2)$.
  Therefore we have 
  \[
    |D_x|^{-\theta+i\kappa} U_\al(t) u_0 
    = \lim_{\del_1,\del_2 \to 0+} \alftr^{-1}( m^{\del_1,\del_2}_t) * u_0.
  \]
  Now we have
  \begin{align*}
    \alftr^{-1}( m^{\del_1,\del_2}_t )(x,y)
    & = c \int_{\R^2} e^{i(x\xi+y\eta)} m^{\del_1,\del_2}_t(\xi,\eta) d\xi d\eta \\
    & = c \int_{\R} |\xi|^{-\theta+i\kappa} e^{-\del_1 \xi^2 + i(x\xi + t\xi|\xi|^\al)}
    \Big( \int_{\R} e^{iy\eta} e^{-(\del_2+\frac{it}\xi)\eta^2} d\eta \Big) d\xi
    \\
    & = c \int_{\R} |\xi|^{-\theta+i\kappa} e^{-\del_1 \xi^2 + ix\xi + it\xi|\xi|^\al}
    (\del_2+\frac{it}\xi)^{-\frac12} e^{-\frac14(\del_2+\frac{it}\xi)^{-1}y^2} d\xi.
  \end{align*}
  Now by the theorem of dominated convergence again we can take the limit
  $\del_2 \to 0+$ in this last expression and get
  \begin{equation*}
    \alftr^{-1}(m_t)
    = c |t|^{-\frac12} \lim_{\del_1\to 0+} 
    \int_{\R} |\xi|^{\frac12-\theta+i\kappa} 
    e^{-\del_1 \xi^2 + i(x\xi + t\xi|\xi|^\al + \sign(\frac{\xi}t)\frac\pi4)}
    e^{\frac{i\xi}{4t} y^2} d\xi.
  \end{equation*}
  Now for $\xi\neq 0$ we set
  $\psi(\xi) := e^{-\del_1 \xi^2 + i \sign(\frac{\xi}t)\frac\pi4}$
  and $\phi(\xi) := \xi|\xi|^\al$.
  Then with our choice of $\theta$ we have
  $|\phi''(\xi)| \sim |\xi|^{1-2\theta}$ and we can use 
  Corollary 2.9 of \cite{kenig:91} to see that
  \begin{equation*}
    |\int_{\R} |\xi|^{\frac12-\theta+i\kappa} 
    e^{-\del_1 \xi^2 + i(x\xi + t\xi|\xi|^\al + \sign(\frac{\xi}t)\frac\pi4)}
    e^{\frac{i\xi}{4t} y^2} d\xi| \ls \jbr{\kappa} |t|^{-\frac12}
  \end{equation*}
  where the implicit constant does not depend on $\delta_1>0$.
  Therefore we get that $\alftr^{-1}(m_t) \in L^\infty(\R^2)$ and 
  \(
    \|\alftr^{-1}(m_t)\|_{L^\infty(\R^2)} \leq C \jbr{\kappa} |t|^{-1}.
  \)
  It follows that we have the decay estimate
  \begin{equation*}
    \| |D_x|^{-\theta+i\kappa} U_\al(t) u_0 \|_{L^\infty} 
    \leq C \jbr{\kappa} |t|^{-1} \|u_0\|_{L^1}
  \end{equation*}
  for all $u_0\in \tf(\R^2)$ and then by continuity also for all
  $u_0\in L^1(\R^2)$.
  By Plancherel we also have that
  $\| |D_x|^{i\kappa} U_\al(t) u_0 \|_{L^2} = \|u_0\|_{L^2}$.
  Now using the interpolation theorem of Stein we get for every
  $2\leq r\leq \infty$ that
  \begin{equation}
    \label{eq:linst_decay}
    \| |D_x|^{(\frac2r-1) \theta} U_\al(t) u_0 \|_{L^r} \leq 
    C |t|^{\frac2r-1} \|u_0\|_{L^{r'}}.
  \end{equation}
  Now \eqref{eq:DGKP_linStrich} follows from \eqref{eq:linst_decay}
  by well-known methods. (See for example \cite{kenig:91}.)
\end{proof}

From this linear version of Strichartz estimates we can deduce the 
following bilinear version
\begin{cor}
  \label{th:DGKP_biStri}
  For $b>\frac12$ we have
  \begin{equation}
    \label{eq:DGKP_biStrich}
    \|u_1 u_2\|_{L^2} \ls \||D_x|^{\frac14-\frac\al 8} u_1\|_{\xsb00b}
                        \||D_x|^{\frac14-\frac\al 8} u_2\|_{\xsb00b}.
  \end{equation}
  Furthermore we have
  \begin{align}
    |\int_{\R^6} \frac{|\xi_1|^{-\frac14+\frac\al 8} |\xi-\xi_1|^{-\frac14+\frac\al 8}}
      {\jbr{\la_1}^b \jbr{\la_2}^b} f_1(\z_1) f_2(\z-\z_1)f_3(\z) d\z_1 d\z|
    & \ls \prod_{i=1}^3 \|f_i\|_{L^2}, \label{eq:DGKP_biStrid} \\
    |\int_{\R^6} \frac{|\xi_1|^{-\frac14+\frac\al 8} |\xi|^{-\frac14+\frac\al 8}}
      {\jbr{\la_1}^b \jbr{\la}^b} f_1(\z_1) f_2(\z-\z_1)f_3(\z) d\z_1 d\z|
    & \ls \prod_{i=1}^3 \|f_i\|_{L^2}, \label{eq:DGKP_dbiStrie} \\
    |\int_{\R^6} \frac{|\xi|^{-\frac14+\frac\al 8} |\xi-\xi_1|^{-\frac14+\frac\al8}}
      {\jbr{\la}^b \jbr{\la_2}^b} f_1(\z_1) f_2(\z-\z_1)f_3(\z) d\z_1 d\z|
    & \ls \prod_{i=1}^3 \|f_i\|_{L^2}. \label{eq:DGKP_dbiStriz} 
  \end{align}
\end{cor}
\begin{proof}
  Setting $r=q=4$ in \eqref{eq:DGKP_linStrich} we get
  \(
    \||D_x|^{-(\frac14-\frac\al8)} U_\al (t)u_0\|_{L^4_{txy}} \ls \|u_0\|_{L^2_{xy}}.
  \)
  Using \cite{ginibre:96}, Lemme 3.3 it follows that
  \(
    \||D_x|^{-(\frac14-\frac\al8)} u\|_{L^4_{txy}} \ls \|u\|_{\xsb00b}
  \)
  or equivalently
  \(
    \|u\|_{L^4_{txy}} \ls \||D_x|^{\frac14-\frac\al8} u\|_{\xsb00b}.
  \)
  Now \eqref{eq:DGKP_biStrich} follows by combining this estimate 
  with H\"older's inequality.
  Setting 
  $f_i(\mu) := |\xi|^{\frac14-\frac{\al}8} \jbr{\la}^b \ftr{u_i}(\mu)$ for
  $i=1,2$ and using duality we see that \eqref{eq:DGKP_biStrich} 
  is equivalent to \eqref{eq:DGKP_biStrid}.
  By suitable changes of variables we also get 
  \eqref{eq:DGKP_dbiStrie} and \eqref{eq:DGKP_dbiStriz}.
\end{proof}

For the part of the product $u_1 u_2$ where the $\xi$-frequency of the 
first factor is significantly smaller than the $\xi$-frequency of the 
second factor we can improve this bilinear Strichartz estimate.
To formulate this improvement let us define for $c>0$
the following operator
\begin{equation}
  \label{eq:DGKP_defparapr}
  \alftr{P_c(u_1,u_2)}(\z) :=
  \int_{\R^3} \chi_{|\xi_1| \leq c|\xi-\xi_1|}(\z,\z_1) 
             \alftr{u_1}(\z_1) \alftr{u_2}(\z-\z_1)\ d\z_1
\end{equation}
We have the following refined bilinear Strichartz estimate
which for the case $\al=2$ was already implicitly used in 
\cite{tzvetkov:99,takaoka:00,tzvetkov:00,takaoka:01_2,takaoka:01,isaza:01}.
\begin{theorem}
  \label{th:DGKP_rebiStri}
  For $b>\frac12$ we have
  \begin{equation}
    \label{eq:DGKP_bilsmooth}
    \|P_{\frac13}(u_1,u_2)\|_{L^2} \ls 
      \||D_x|^{\frac12} u_1\|_{\xsb00b}
      \||D_x|^{-\frac\al 4} u_2\|_{\xsb00b}
  \end{equation}
\end{theorem}
\
For the proof of the theorem we need the following Lemma
\begin{lemma}
  \label{le:DGKP_absres}
  For $\al>0$ set $\phi_\al(\xi) := \xi |\xi|^\al$ and
  \begin{equation}
    \label{eq:DGKP_defres}
    r_\al(\xi,\xi_1) := \phi_\al(\xi) - \phi_\al(\xi_1) - \phi_\al(\xi-\xi_1),
    \quad \xi,\xi_1\in\R. 
  \end{equation}
  We then have for every $\xi,\xi_1\in\R$
  \begin{equation}
    \label{eq:DGKP_absres}
    \frac\al{2^\al} |\xi_{\min}| |\xi_{\max}|^\al \leq
    |r_\al(\xi,\xi_1)| \leq (\al+1+\frac1{2^\al}) |\xi_{\min}| |\xi_{\max}|^\al.
  \end{equation}
\end{lemma}
\
\begin{proof}[Proof of Lemma~\ref{le:DGKP_absres}]
  Suppose first that $|\xi_{\min}|=|\xi_1|$.
  Then we have
  \[
  |\phi_\al(\xi_1)|=|\xi_{\min}|^{\al+1} \leq \frac1{2^\al}|\xi_{\min}||\xi_{\max}|^\al
  \]
  because $|\xi_{\min}|\leq \frac12 |\xi_{\max}|$ and
  \[
  |\phi_\al(\xi) - \phi_\al(\xi-\xi_1)| 
  = |\phi'_\al(\xi-\theta\xi_1)||\xi_1|
  = (\al+1) |\xi-\theta\xi_1|^\al |\xi_{\min}|
  \]
  for some $\theta\in [0,1]$.
  As we have $|\xi_1| \leq |\xi|$ it follows that
  \[
  \min_{\theta\in [0,1]} |\xi-\theta \xi_1| 
  = \min \{|\xi|,|\xi-\xi_1|\} = |\xi_{\med}| 
  \geq \frac12 |\xi_{\max}|
  \]
  and 
  \[
  \max_{\theta\in [0,1]} |\xi-\theta \xi_1| 
  = \max \{|\xi|,|\xi-\xi_1|\} = |\xi_{\max}|.
  \]
  Putting these estimates together we get
  \begin{align*}
    |r_\al(\xi,\xi_1)| 
    & \geq |\phi_\al(\xi) - \phi_\al(\xi-\xi_1)| - |\phi_\al(\xi_1)| \\
    & \geq (\al+1)\frac1{2^\al}|\xi_{\max}|^\al|\xi_{\min}| -
    \frac1{2^\al}|\xi_{\min}||\xi_{\max}|^\al\\
    & = \frac\al{2^{\al}}|\xi_{\min}||\xi_{\max}|^\al
  \end{align*}
  and
  \begin{align*}
    |r_\al(\xi,\xi_1)| 
    & \leq |\phi_\al(\xi) - \phi_\al(\xi-\xi_1)| + |\phi_\al(\xi_1)| \\
    & \leq (\al+1)|\xi_{\max}|^\al|\xi_{\min}| + \frac1{2^\al}|\xi_{\min}||\xi_{\max}|^\al\\
    & = (\al+1+\frac1{2^{\al}})|\xi_{\min}||\xi_{\max}|^\al
  \end{align*}
  which proves \eqref{eq:DGKP_absres} in the case $|\xi_{\min}|=|\xi_1|$.
  Noting that
  $r_\al(\xi,\xi_1)=r_\al(\xi,\xi-\xi_1)=-r_\al(\xi-\xi_1,\xi)$
  we see that we also get \eqref{eq:DGKP_absres} in the other cases.
\end{proof}

\begin{proof}[Proof of Theorem~\ref{th:DGKP_rebiStri}]
  Let $f_1(\z) := |\xi|^{\frac12} \jbr{\la}^b \ftr{u_1}(\z)$
  and $f_2(\z) := |\xi|^{-\frac{\al}4} \jbr{\la}^b \ftr{u_2}(\z)$.
  We have to show that
  \begin{equation*}
    \| \int_{\R^3} \chi_{|\xi_1| \leq \frac13|\xi-\xi_1|} 
      \frac{|\xi_1|^{-\frac12} |\xi-\xi_1|^{\frac{\al}4}}
           {\jbr{\la_1}^b {\jbr{\la_2}^b}}
      f_1(\z_1) f_2(\z-\z_1) d\z_1 \|_{L^2_\z} 
      \ls \|f_1\|_{L^2} \|f_2\|_{L^2}
  \end{equation*}
  which by duality is equivalent to
  \begin{equation}
    \label{eq:DGKP_bilsmdual}
    | \int_{\R^6}\ \chi_{|\xi_1| \leq \frac13|\xi-\xi_1|} 
    \frac{|\xi_1|^{-\frac12} |\xi-\xi_1|^{\frac{\al}4}}
    {\jbr{\la_1}^b {\jbr{\la_2}^b}} f_1(\z_1) f_2(\z-\z_1)f_3(\z) d\z_1 d\z |
    \ls \prod_{i=1}^3 \|f_i\|_{L^2}.
  \end{equation}
  By use of the Cauchy-Schwarz inequality it suffices to show that 
  $\sup_\z I(\z)^{\frac12} < \infty$ where
  \begin{equation*}
    I(\z) := \int_{\R^3} \chi_{|\xi_1| \leq \frac13|\xi-\xi_1|} 
      \frac{|\xi_1|^{-1} |\xi-\xi_1|^{\frac{\al}2}}
           {\jbr{\la_1}^{2b} {\jbr{\la_2}^{2b}}}\ d\z_1.
  \end{equation*}
  For fixed $\z$ we now use the change of variables 
  $T : \z_1 \mapsto (\nu,\la_1,\la_2)$
  where
  \[ \nu(\z_1) := r_\al(\xi,\xi_1)
  = \xi |\xi|^\al - \xi_1 |\xi_1|^\al - (\xi-\xi_1) |\xi-\xi_1|^\al.
  \]
  Let us also recall the definition of $\la_1$ and $\la_2$
  \begin{align*}
    \la_1(\z_1) & = \tau_1 - \xi_1 |\xi_1|^\al + \frac{\eta_1^2}{\xi_1},\\
    \la_2(\z_1) & = \tau-\tau_1 - (\xi-\xi_1)|\xi-\xi_1|^\al 
                 + \frac{(\eta-\eta_1)^2}{\xi-\xi_1}.
  \end{align*}
  Observe that 
  \begin{equation}
    \label{eq:DGKP_resonid}
    \la_1 + \la_2 - \la 
    = \nu + \frac{(\xi\eta_1 - \eta\xi_1)^2}{\xi\xi_1(\xi-\xi_1)}.
  \end{equation}
  Therefore we have 
  \begin{equation*}
    |\partial_{\eta_1} (\la_1 + \la_2)| = 
    2|\xi \frac{\xi\eta_1-\eta\xi_1}{\xi\xi_1(\xi-\xi_1)}|
    = 2 \frac{|\xi|^{\frac12} |\la_1+\la_2-\la-\nu|^{\frac12}}
    {|\xi_1|^{\frac12} |\xi-\xi_1|^{\frac12}}.
  \end{equation*}
  Furthermore we have
  \(
    \partial_{\xi_1} \nu = (\al+1)(|\xi-\xi_1|^\al - |\xi_1|^\al). 
  \)
  As we only consider the region where $|\xi_1|\leq \frac13|\xi-\xi_1|$, 
  i. e. $|\xi_1| = |\xi_{\min}|$ and $|\xi-\xi_1| \sim |\xi_{\max}|$ 
  we have by \eqref{eq:DGKP_absres} that 
  $|\nu| \sim |\xi_1||\xi-\xi_1|^\al$.
  We also have $|\partial_{\xi_1} \nu| \gs |\xi-\xi_1|^\al$
  in this region.
  Therefore we have that
  \begin{align*}
    |D_{\z_1} T| 
    & = |\partial_{\xi_1} \nu|
        |\partial_{\eta_1} \la_1+\partial_{\eta_1} \la_2| 
      \gs |\xi_1|^{-\frac12} |\xi-\xi_1|^\al |\la_1+\la_2-\la-\nu|^{\frac12} \\
    & \gs |\xi_1|^{-1} |\xi-\xi_1|^{\frac\al 2} |\nu|^{\frac12}
          |\la_1+\la_2-\la-\nu|^{\frac12}.
  \end{align*}
  Let us notice that it is possible to divide the region of integration  
  into a finite number of open subsets $U_i$ such that $T$ is 
  an injective $C^1$-function in $U_i$ with nonvanishing Jacobian.
  As we are in the KP-II-case both terms on the right hand side of
  \eqref{eq:DGKP_resonid} have the same sign which implies that
  $|\nu| \leq |\la_1+\la_2-\la|$.
  So performing the change of variables and using the following
  elementary inequality
  \begin{equation*}
    \int_{-K}^K \frac{d\nu}{|\nu|^\frac12 |a-\nu|^\frac12}
    \ls \frac{K^\frac12}{|a|^\frac12},\quad a\neq 0
  \end{equation*}
  we get
  \begin{equation*}
    I(\z)\ls \int_{\R^3}
    \frac{\chi_{|\nu|\leq |\la_1+\la_2-\la|} d\nu d\la_1 d\la_2}
         {\jbr{\la_1}^{2b}\jbr{\la_2}^{2b}|\nu|^{\frac12}|\la_1+\la_2-\la-\nu|^{\frac12}}
    \ls \int_{\R^2} \frac{d\la_1 d\la_2}{\jbr{\la_1}^{2b}\jbr{\la_2}^{2b}}
    \ls 1.
  \end{equation*}
\end{proof}

\begin{rem}
  In fact we get $\eqref{eq:DGKP_bilsmdual}$ also without the cut-off function
  $\chi_{|\xi_1| \leq \frac13 |\xi-\xi_1|}$,
  as in the region where $|\xi_1| > \frac13 |\xi-\xi_1|$ we have that 
  \(
    |\xi_1|^{-\frac12} |\xi-\xi_1|^{\frac\al 4}
    \ls |\xi_1|^{-\frac14+\frac\al 8} |\xi-\xi_1|^{-\frac14+\frac\al 8}
  \)
  and so the estimate in this region follows from the bilinear 
  Strichartz estimate \eqref{eq:DGKP_biStrid}.
  Like in the case of the bilinear Strichartz estimates we also get the
  following dual versions of estimate \eqref{eq:DGKP_bilsmdual} 
  (without the cut-off function) by an appropriate change of variables
  \begin{align}
    | \int_{\R^6} \frac{|\xi_1|^{-\frac12} |\xi|^{\frac{\al}4}}
      {\jbr{\la_1}^b {\jbr{\la}^b}} f_1(\z_1) f_2(\z-\z_1)f_3(\z) d\z_1 d\z | 
      & \ls \prod_{i=1}^3 \|f_i\|_{L^2}, \label{eq:DGKP_dbilsmd} \\
    | \int_{\R^6} \frac{|\xi|^{-\frac12} |\xi-\xi_1|^{\frac{\al}4}}
      {\jbr{\la}^b {\jbr{\la_2}^b}} f_1(\z_1) f_2(\z-\z_1)f_3(\z) d\z_1 d\z | 
      & \ls \prod_{i=1}^3 \|f_i\|_{L^2}, \label{eq:DGKP_dbilsmdz} \\
    | \int_{\R^6} \frac{|\xi|^{-\frac12} |\xi_1|^{\frac{\al}4}}
      {\jbr{\la}^b {\jbr{\la_1}^b}} f_1(\z_1) f_2(\z-\z_1)f_3(\z) d\z_1 d\z | 
      & \ls \prod_{i=1}^3 \|f_i\|_{L^2}. \label{eq:DGKP_dbilsmdd}    
  \end{align}
\end{rem}

\section{The main bilinear estimate}
\label{se:DGKP_mbe}
\
In the following formulation and proof of the crucial
bilinear estimate needed to prove Theorem~\ref{th:DGKP_maindgkp}
we will only consider the case $s_2=0$ (and write $s$ for $s_1$)
to simplify the presentation.
Note that the case $s_2>0$ follows from this special case, as 
in the general case we only get an extra term 
\(
  \frac{\jbr{\eta}^{s_2}}{\jbr{\eta_1}^{s_2}\jbr{\eta-\eta_1}^{s_2}}
\)
in the integral inequalities we have to prove 
(see \eqref{eq:DGKP_intesttp}).
But this term is always bounded above for $s_2\geq 0$.
\begin{theorem}
  \label{th:DGKP_bilest}
  For $\frac43<\al\leq 6$ and $s>\max(1-\frac34 \al,\frac14-\frac38 \al)$ 
  let us choose $0<\eps\leq \frac1{16}$ such that
  \[
  \eps \leq \min(\frac1{\al+1}(s+\frac34 \al - 1),
                 \frac1{\al+1}(\frac34 \al - 1),
                 \frac1\al(2s+\frac34 \al - \frac12))
  \]
  and set $b:=\frac12 +\frac\eps2$, $b'=-\frac12 + \eps$,
  $b_1:=\max(\frac3{2\al}-\frac34+\eps,0)$ and
  $\lf=\frac12+b_1$.
  Let us define $X=X_1+X_2$ and $\tilde{X}=\tilde{X}_1+\tilde{X}_2$ where
  \begin{align*}
    X_1 := \xsb{s}0{b-b'},\quad & 
    X_2 := \xsb[\lf]s0b \cap \xsb[\lf]{s-(\al+1)b_1}0{b+b_1}, \\
    \tilde{X}_1 := \xsb{s}00, \quad &
    \tilde{X}_2 := \xsb[\lf]s0{b'} \cap \xsb[\lf]{s-(\al+1)b_1}0{b'+b_1}
  \end{align*}
  We then have
  \begin{equation}
    \label{eq:DGKP_nonlinest}
    \|\partial_x (u_1 u_2)\|_{\tilde{X}} \leq C \|u_1\|_{X} \|u_2\|_{X}.
  \end{equation}
\end{theorem}

\begin{rem}
  The spaces $X$ and $\tilde{X}$ defined in Theorem~\ref{th:DGKP_bilest} 
  are built by taking sums and intersections of the Bourgain type
  spaces of Section~\ref{se:DGKP_defxsb}.
  Therefore it is easy to see that they also satisfy the linear 
  estimates of Propositions~\ref{pr:DGKP_linesth} and 
  \ref{pr:DGKP_linesti}, i. e. 
  \begin{equation}
    \label{eq:DGKP_linestivX}
    \| \psi U_\al (t) u_0 \|_X \ls \|u_0\|_{H^{s,0}(\R^2)}
  \end{equation}
  and 
  \begin{equation}
    \label{eq:DGKP_linestiX}
    \| \psi_T \int_0^t U(t - t') F(t') dt' \|_X
    \ls T^{1-(b-b')} \|F\|_{\tilde{X}}.
  \end{equation}  
\end{rem}

\begin{rem}
  The sum structure of the spaces $X$ and $\tilde{X}$ is the essential
  ingredient to use the additional weight 
  $(\frac{\jbr{\xi}}{|\xi|})^\lf$, 
  which is incorporated in the definition of $X_2$ and $\tilde{X}_2$ (see
  \eqref{eq:DGKP_xsbndef}), to lower the $x$-regularity $s$ in the 
  bilinear estimate without imposing a low frequency condition on the
  initial data.
  Therefore, in the case $\al=2$ of the Kadomtsev-Petviashvili II equation
  we are able to show the local well-posedness for all
  $s>-\frac12$ without a low frequency condition on the initial
  data whereas the counterexamples in \cite{takaoka:01}
  show that it is not possible to get the bilinear estimate for 
  $-\frac12<s<-\frac13$ and $\sigma=0$.
\end{rem}

\begin{rem}
  \label{rm:DGKP_xtsp}
  In the case $\al=4$ of the fifth order KP II equation it is
  possible to get the bilinear estimate \eqref{eq:DGKP_nonlinest} 
  in the spaces $X=\xsb{s}{0}{b}$ and $\tilde{X}=\xsb{s}{0}{b'}$ 
  (i. e. choosing $b_1=0$ and $\sigma=0$) as is shown in \cite{isaza:06}.
  More generally, this is true for all $\al>\frac52$ which can be
  seen by refining the estimate of Lemma~\ref{le:DGKP_bilinqee}
  by an additional dyadic decomposition and interpolation
  argument as used in \cite{takaoka:01}, pp. 89-92.
\end{rem}

To prove Theorem~\ref{th:DGKP_bilest}, we will split the nonlinear term 
$\partial_x (u_1 u_2)$ into various pieces and give estimates
in appropriate $X^{s,b}_\sigma$-spaces for each of these pieces 
(see Lemmas~\ref{le:DGKP_bilinqnn}-\ref{le:DGKP_bilinqzz}).
We will then combine these estimates to prove \eqref{eq:DGKP_nonlinest}.
First of all, with $P_c$ defined as in \eqref{eq:DGKP_defparapr},
we can write
\begin{equation*}
  \partial_x (u_1 u_2) 
  = \partial_x P_1(u_1,u_2) + \partial_x P_1(u_2,u_1)
\end{equation*}
As the main bilinear estimate \eqref{eq:DGKP_nonlinest} is symmetric in 
$u_1$ and $u_2$, it suffices to prove it only for $\partial_x P_1(u_1,u_2)$.
This expression can be decomposed further into
\begin{equation}
  \label{eq:DGKP_nonldecmp}
  \partial_x P_1(u_1 u_2)
  = Q_{00}(u_1,u_2) + \sum_{i=1}^2 \sum_{j=0}^2 Q_{ij}(u_1,u_2).
\end{equation}
The operators $Q_{ij}$ are defined by
\begin{equation*}
  \alftr{Q_{ij}(u_1,u_2)}(\z) 
  = i\xi \int_{\R^3} \chi_{A_{ij}}(\z_1,\z) 
           \alftr{u_1}(\z_1) \alftr{u_2}(\z-\z_1) d\z_1.
\end{equation*}
where $A_{00} := \{ (\z_1,\z)\in \R^6 \mid |\xi_1|\leq |\xi-\xi_1|\leq 1 \}$
and $A_{ij} := \Xi_i \cap \Lambda_j$ for $1\leq i\leq 2$, $0\leq j\leq 2$ with
\begin{align*}
  \Xi_1 & = 
  \{ (\z_1,\z)\in \R^6 \mid |\xi_1| \leq \frac13 |\xi-\xi_1|, |\xi-\xi_1|\geq 1 \}, \\
  \Xi_2 & = \{ (\z_1,\z)\in \R^6 \mid \frac13 |\xi-\xi_1|\leq |\xi_1| \leq |\xi-\xi_1|,
          |\xi-\xi_1|\geq 1 \}, \\
  \Lambda_0 & = \{ (\z_1,\z)\in \R^6 \mid |\la| = |\la_{\max}| \} \\
  \Lambda_j & = \{ (\z_1,\z)\in \R^6 \mid |\la_j| = |\la_{\max}| \}\quad (j=1,2)
\end{align*}
Let us explain what the meaning of the regions $\Xi_1$ and $\Xi_2$
is.
In $\Xi_1$ we have that $2\leq 2|\xi-\xi_1|\leq 3|\xi|\leq 4|\xi-\xi_1|$,
i. e. $\xi$ and $\xi-\xi_1$ are comparable in size and are both
bounded away from zero, whereas $\xi_1$ is the smallest of the 
frequencies dual to the $x$-variable.
In $\Xi_2$ we have that $\xi_1$ and $\xi-\xi_1$ are comparable in size
and are both bounded away from zero, whereas $\xi$ may be small here.
For each of the operators $Q_{ij}$ we will now show estimates of the form
\begin{equation*}
  \|Q_{ij}(u_1,u_2)\|_{\xsb[\lf']{s'}0{b'}} 
  \ls \|u_1\|_{\xsb[\lf_1]{s_1}0{b_1}} \|u_2\|_{\xsb[\lf_2]{s_2}0{b_2}}.
\end{equation*}
By definition \eqref{eq:DGKP_xsbndef} of the $X^{s,b}_\lf$-norm
and setting 
$f_l(\z) := |\xi|^{-\lf} \jbr{\xi}^{s_l+\lf} \jbr{\la}^{b_l} \ftr{u_l}(\z)$ 
this is equivalent to
\begin{equation*}
  \left\|
    \frac{|\xi|\jbr{\xi}^{s'+\lf'}}{|\xi|^{\lf'}\jbr{\la}^{-b'}}  
    \int\limits_{\R^3} \chi_{A_{ij}}(\mu_1,\mu)
    \frac{|\xi_1|^{\lf_1} f_1(\z_1) |\xi-\xi_1|^{\lf_2} f_2(\z-\z_1) d\z_1}
    {\jbr{\xi_1}^{s_1+\lf_1} \jbr{\la_1}^{b_1} \jbr{\xi-\xi_1}^{s_2+\lf_2}
    \jbr{\la_2}^{b_2}}
  \right\|_{L^2_\z} \ls \|f_1\|_{L^2} \|f_2\|_{L^2}. 
\end{equation*}
Using duality this estimate is equivalent to
\begin{equation}
  \label{eq:DGKP_intesttp}
    \left|
      \int_{A_{ij}} 
      \frac{|\xi|\jbr{\xi}^{s'+\lf'} |\xi_1|^{\lf_1} |\xi-\xi_1|^{\lf_2} 
             f_1(\z_1)f_2(\z-\z_1)f_3(\z)}
           {|\xi|^{\lf'} \jbr{\xi_1}^{s_1+\lf_1} \jbr{\xi-\xi_1}^{s_2+\lf_2}
             \jbr{\la}^{-b'} \jbr{\la_1}^b \jbr{\la_2}^b}
      d\z d\z_1
    \right| \ls \prod_{i=1}^3 \|f_i\|_{L^2}.  
\end{equation}
The main ingredients we use in the proof of these estimates are 
the bilinear Strichartz estimates
of corollary~\ref{th:DGKP_biStri} and Theorem~\ref{th:DGKP_rebiStri}
and a use of the ``resonance identity'' \eqref{eq:DGKP_resonid}.
We already noted that the two terms on the right hand side of
\eqref{eq:DGKP_resonid} have the same sign. 
Therefore we have
\begin{equation}
  \label{eq:DGKP_lmeste}
  |\la_{\max}| \geq \frac13 |\la_1+\la_2-\la| 
  \geq \frac13 |\nu|
  \geq \frac\al{3\cdot 2^\al} |\xi_{\min}| |\xi_{\max}|^\al.
\end{equation}
where for the last inequality we used \eqref{eq:DGKP_absres}.

\begin{lemma} 
  \label{le:DGKP_bilinqnn}
  We have that
  \begin{equation}
    \label{eq:DGKP_bilinqnn}
    \|Q_{00}(u_1,u_2)\|_{\xsb{s}00} 
    \ls \|u_1\|_{\xsb{s}0b} \|u_2\|_{\xsb{s}0b}
  \end{equation}
  provided that $b>\frac12$, $\al\leq 6$ and $s\in\R$.
\end{lemma}
\begin{proof}
  We have to prove that
  \begin{equation*}
    |\int_{A_{00}} k_{00}(\z_1,\z)
    \frac{|\xi_1|^{-\frac12} |\xi-\xi_1|^{\frac\al4}}{\jbr{\la_1}^b \jbr{\la_2}^b}
    f_1(\z_1)f_2(\z-\z_1)f_3(\z) d\z d\z_1|
    \ls \prod_{i=1}^3 \|f_i\|_{L^2}
  \end{equation*}
  where 
  \(
  k_{00}(\z_1,\z)=\jbr{\xi}^s \jbr{\xi_1}^{-s} \jbr{\xi-\xi_1}^{-s}
  |\xi| |\xi_1|^{\frac12} |\xi-\xi_1|^{-\frac\al4}.
  \)
  On $A_{00}$ we have that $|\xi_1|\leq |\xi-\xi_1|\leq 1$ and therefore also
  $|\xi|\leq 2|\xi-\xi_1|\leq 2$, so that
  \(
    k_{00}(\z_1,\z) \ls |\xi| |\xi_1|^{\frac12} |\xi-\xi_1|^{-\frac\al4}
    \ls |\xi-\xi_1|^{\frac32-\frac\al 4} \ls 1.
  \)
  where the last inequality follows from $\al\leq 6$.
  Therefore \eqref{eq:DGKP_bilinqnn} follows from the refined
  bilinear Strichartz estimate \eqref{eq:DGKP_bilsmdual}.
\end{proof}

\begin{lemma} 
  \label{le:DGKP_bilinqen}
  We have that
  \begin{equation}
    \label{eq:DGKP_bilinqen}
    \|Q_{10}(u_1,u_2)\|_{\xsb[1]{s-(\al+1)b_1}0{b'+b_1}} 
    \ls \|u_1\|_{\xsb{s}0b} \|u_2\|_{\xsb{s}0b}
  \end{equation}
  provided that $b>\frac12$, $b'>-\frac12$ and
  \begin{gather}
    0 \leq b_1 \leq -b', \label{eq:DGKP_qebe} \\
    b'\leq \frac1{\al+1}(\min(0,s)-\frac32+\frac\al4). \label{eq:DGKP_qebz}
  \end{gather}
\end{lemma}
\begin{proof}
  We have to prove that
  \begin{equation*}
    |\int_{A_{10}} k_{10}(\z_1,\z) \frac{|\xi_1|^{-\frac12} |\xi-\xi_1|^{\frac\al 4}}
                     {\jbr{\la_1}^b \jbr{\la_2}^b}
                f_1(\z_1)f_2(\z-\z_1)f_3(\z) d\z d\z_1|
    \ls \prod_{i=1}^3 \|f_i\|_{L^2}
  \end{equation*}
  where
  \(
    k_{10}(\z_1,\z) 
    = \jbr{\la}^{b'+b_1} 
      \jbr{\xi}^{1+s-(\al+1)b_1} \jbr{\xi_1}^{-s} \jbr{\xi-\xi_1}^{-s}
      |\xi_1|^{\frac12} |\xi-\xi_1|^{-\frac\al 4}.
  \)
  We show that $k_{10}$ is bounded in $A_{10}$, then the lemma follows
  from the refined bilinear Strichartz estimate \eqref{eq:DGKP_bilsmdual}.
  In region $A_{10}$ we have 
  $1\leq |\xi-\xi_1|\sim |\xi| \sim \jbr{\xi}$
  and $|\la|=|\la_{\max}|\gs |\xi_1||\xi|^\al$,
  so using \eqref{eq:DGKP_qebe} we get 
  \( 
    k_{10}(\z_1,\z)
    \ls |\xi|^{1 -\frac\al 4+\al b'-b_1} \jbr{\xi_1}^{-s} |\xi_1|^{\frac12+b'+b_1}
  \)
  Because of \eqref{eq:DGKP_qebe} we have $\frac12+b'+b_1\geq 0$.
  So using $|\xi_1|\leq |\xi| \sim \jbr{\xi}$ it follows that
  \(
    k_{10}(\z_1,\z)
    \ls \jbr{\xi}^{\frac32-\frac\al4+(\al+1)b'-\min(0,s)} \ls 1
  \)
  where the last inequality follows from \eqref{eq:DGKP_qebz}.
\end{proof}

\begin{lemma}
  \label{le:DGKP_bilinqez}
  We have that
  \begin{equation}
    \label{eq:DGKP_bilinqez}
    \|Q_{12}(u_1,u_2)\|_{\xsb[1]{s-(\al+1)b_1}0{b'+b_1}} 
    \ls \|u_1\|_{\xsb{s}0b} \|u_2\|_{\xsb{s}0b}
  \end{equation}
  provided that $b>\frac12$, $b'>-\frac12$ and \eqref{eq:DGKP_qebe} and
  \eqref{eq:DGKP_qebz} hold.
\end{lemma}
\begin{proof}
  We have to show that
  \begin{equation*}
    |\int_{A_{12}} k_{12}(\z_1,\z) \frac{|\xi_1|^{-\frac12} |\xi|^{\frac\al 4}}
                     {\jbr{\la_1}^b \jbr{\la}^b}
                f_1(\z_1)f_2(\z-\z_1)f_3(\z) d\z d\z_1|
    \ls \prod_{i=1}^3 \|f_i\|_{L^2}
  \end{equation*}
  where
  \(
    k_{12}(\z_1,\z) 
    = \jbr{\la}^{b'+b_1+b} \jbr{\la_2}^{-b}
      \jbr{\xi}^{1+s-(\al+1)b_1} \jbr{\xi_1}^{-s} \jbr{\xi-\xi_1}^{-s}
      |\xi_1|^{\frac12} |\xi|^{-\frac\al 4}.
  \)  
  Now using $b'+b_1+b\geq 0$ and $|\la|\leq |\la_2|$ in $A_{12}$
  we get $\jbr{\la}^{b'+b_1+b} \jbr{\la_2}^{-b} \leq \jbr{\la_2}^{b'+b_1}$.
  Now the boundedness of $k_{12}$ on $A_{12}$ follows exactly like the
  boundedness of $k_{10}$ in Lemma~\ref{le:DGKP_bilinqen}.
  Then \eqref{eq:DGKP_bilinqez} follows from the refined bilinear 
  Strichartz estimate \eqref{eq:DGKP_dbilsmd}.
\end{proof}

\begin{lemma}
  \label{le:DGKP_bilinqee}
  We have that
  \begin{equation}
    \label{eq:DGKP_bilinqeea}
    \|Q_{11}(u_1,u_2)\|_{\xsb[1]s00} 
    \ls \|u_1\|_{\xsb{s}0{b-b'}} \|u_2\|_{\xsb{s}0b}    
  \end{equation}
  and
  \begin{equation}
    \label{eq:DGKP_bilinqeeb}
    \|Q_{11}(u_1,u_2)\|_{\xsb[1]{s-(\al+1)b_2}0{b'+b_2}} 
    \ls \|u_1\|_{\xsb[\lf]{s-(\al+1)b_1}0{b+b_1}} \|u_2\|_{\xsb{s}0b}
  \end{equation}
  provided that $b>\frac12$, $b'>-\frac12$, \eqref{eq:DGKP_qebe} and
  \eqref{eq:DGKP_qebz} hold and 
  \begin{gather}
    0\leq b_2\leq b_1, \label{eq:DGKP_qebd} \\
    \lf \geq b_1-b' \geq \frac3{2\al} - \frac14. \label{eq:DGKP_qebv}  
  \end{gather}
\end{lemma}
\begin{proof}
  We first show \eqref{eq:DGKP_bilinqeea}.
  We have to prove that
  \begin{equation*}
    |\int_{A_{11}} \tilde{k}_{11}(\z_1,\z) 
    \frac{|\xi_1|^{-\frac12} |\xi-\xi_1|^{\frac\al 4}}
         {\jbr{\la_1}^b \jbr{\la_2}^b} 
    f_1(\z_1)f_2(\z-\z_1)f_3(\z) d\z d\z_1|
    \ls \prod_{i=1}^3 \|f_i\|_{L^2}
  \end{equation*}
  where 
  \(
    \tilde{k}_{11}(\z_1,\z)
    =\jbr{\la_1}^{b'}   
    \jbr{\xi}^{1+s} \jbr{\xi_1}^{-s} \jbr{\xi-\xi_1}^{-s}
    |\xi_1|^{\frac12} |\xi-\xi_1|^{-\frac\al 4}.
  \)
  We will show that $\tilde{k}_{11}$ is bounded, then 
  \eqref{eq:DGKP_bilinqeea} follows from the refined bilinear
  Strichartz estimate \eqref{eq:DGKP_bilsmdual}.
  In the region $A_{11}$ we have $1\leq |\xi-\xi_1| \sim |\xi| \sim \jbr{\xi}$
  and $|\la_1|=|\la_{\max}|\gs |\xi_1| |\xi|^\al$.
  Because of $b'\leq 0$ we have
  $\tilde{k}_{11}(\z_1,\z)
  \ls |\xi|^{1-\frac\al 4+\al b'} \jbr{\xi_1}^{-s} |\xi_1|^{\frac12+b'}\ls 1$,
  where the last inequality follows exactly as in the proof of 
  Lemma~\ref{le:DGKP_bilinqen}.

  Now we show \eqref{eq:DGKP_bilinqeeb}.
  We have to prove that
  \begin{equation*}
    |\int_{A_{11}} k_{11}(\z_1,\z)
    \frac{|\xi|^{-\frac14+\frac\al 8} |\xi-\xi_1|^{-\frac14+\frac\al 8}}
         {\jbr{\la}^b \jbr{\la_2}^b}
    f_1(\z_1)f_2(\z-\z_1)f_3(\z) d\z d\z_1| \ls \prod_{i=1}^3 \|f_i\|_{L^2}
  \end{equation*}
  where
  \begin{equation*}
    k_{11}(\z_1,\z) = \frac{
    \jbr{\xi}^{1+s-(\al+1)b_2} \jbr{\xi_1}^{-s+(\al+1)b_1-\lf} \jbr{\xi-\xi_1}^{-s}
    |\xi|^{\frac14-\frac\al 8} |\xi_1|^\lf |\xi-\xi_1|^{\frac14-\frac\al 8}}
    {\jbr{\la}^{-b'-b_2-b} \jbr{\la_1}^{b+b_1}}.
  \end{equation*}
  Now \eqref{eq:DGKP_bilinqeeb} follows from the bilinear Strichartz
  estimate \eqref{eq:DGKP_biStrid}
  if we show that $k_{11}$ is bounded on $A_{11}$.
  Because of $|\la|\leq |\la_1|$ and $-b'-b_2-b \leq 0$, we have
  \[
  \jbr{\la}^{-b'-b_2-b} \jbr{\la_1}^{b+b_1} \geq \jbr{\la_1}^{-b'+b_1-b_2}.
  \]
  In $A_{11}$ we have $1\leq |\xi-\xi_1| \sim |\xi| \sim \jbr{\xi}$ and
  $|\la_1| \gs |\xi_1| |\xi|^\al$, so using $-b'+b_1-b_2 \geq 0$ we get
  \begin{equation*}
    k_{11}(\z_1,\z) \ls
     |\xi|^{\frac32-\frac\al 4+\al(b'-b_1)-b_2} 
     \jbr{\xi_1}^{-s+(\al+1)b_1-\lf} |\xi_1|^{\lf +b'-b_1+b_2}.
  \end{equation*}
  Because of $|\xi_1|\leq |\xi|$ in $A_{11}$, we have 
  $|\xi|^{-b_2}|\xi_1|^{b_2}\leq 1$.
  By \eqref{eq:DGKP_qebv} we have $\lf + b' - b_1 \geq 0$.
  If $|\xi_1|\leq 1$ we have 
  \(
    k_{11}(\z_1,\z)\ls |\xi|^{ \frac32 - \frac\al 4 + \al (b'-b_1)} \ls 1
  \)
  where the last inequality follows by \eqref{eq:DGKP_qebv}.
  If $|\xi_1|\geq 1$ we have $\jbr{\xi_1} \sim |\xi_1|$ and therefore
  \begin{equation*}
    k_{11}(\z_1,\z)
    \ls |\xi|^{\frac32-\frac\al 4+\al(b'-b_1)} \jbr{\xi_1}^{-s + b' + \al b_1}
  \end{equation*}
  Now if $-s + b' + \al b_1 \leq 0$ this term is bounded
  because of \eqref{eq:DGKP_qebv} as above.
  If $-s + b' + \al b_1 > 0$ this term is bounded by 
  $c |\xi|^{\frac32-\frac\al 4-s+(\al+1)b'}$ which is bounded 
  because of \eqref{eq:DGKP_qebz}.
\end{proof}

\begin{lemma} 
  \label{le:DGKP_bilinqzn}
  We have that
  \begin{equation}
    \label{eq:DGKP_bilinqzn}
    \|Q_{20}(u_1,u_2)\|_{\xsb[\lf]{s-(\al+1)b_2}0{b'+b_2}} 
    \ls \|u_1\|_{\xsb{s}0b} \|u_2\|_{\xsb{s}0b}
  \end{equation}
  provided that $b>\frac12$, $b'>-\frac12$, \eqref{eq:DGKP_qebe} and 
  \eqref{eq:DGKP_qebz} hold and
  \begin{gather}
    0\leq b_2\leq b_1\leq -b'+ \frac1\al (2s-\frac12+\frac\al4), 
    \label{eq:DGKP_qzbe} \\
    \lf \leq 1 + b' + b_1. \label{eq:DGKP_qzbz}
  \end{gather}
\end{lemma}
\begin{proof}
  We have to show that
  \begin{equation*}
    |\int_{A_{20}} k_{20}(\z_1,\z) 
      \frac{|\xi_1|^{-\frac14+\frac\al 8} |\xi-\xi_1|^{-\frac14+\frac\al 8}}
           {\jbr{\la_1}^b \jbr{\la_2}^b}
      f_1(\z_1)f_2(\z-\z_1)f_3(\z) d\z d\z_1|
    \ls \prod_{i=1}^3 \|f_i\|_{L^2}
  \end{equation*}
  where
  \(
    k_{20}(\z_1,\z) = \jbr{\la}^{b'+b_2}
    \jbr{\xi}^{s-(\al+1)b_2+\lf} \jbr{\xi_1}^{-s} \jbr{\xi-\xi_1}^{-s} 
    |\xi|^{1-\lf} |\xi_1|^{\frac14-\frac\al 8} |\xi-\xi_1|^{\frac14-\frac\al 8}.
  \)
  Now if we show that $k_{20}$ is bounded on $A_{20}$, the lemma follows
  from \eqref{eq:DGKP_biStrid}.
  In $A_{20}$ we have $1\leq |\xi-\xi_1|\sim |\xi_1| \sim \jbr{\xi_1}$
  and $|\la|=|\la_{\max}|$, so
  \begin{equation}
    \label{eq:DGKP_exprxiz}
    k_{20}(\z_1,\z) \ls h(\z_1,\z) :=
      \jbr{\la_{\max}}^{b'+b_2} \jbr{\xi}^{s-(\al+1)b_2+\lf}
      |\xi|^{1-\lf} |\xi_1|^{-2s+\frac12-\frac\al4}.
  \end{equation}
  We will now show that $h$ is bounded in $\Xi_2$.
  Let us first consider the case that $|\xi|\geq 1$.
  Then because of $|\la_{\max}| \geq |\xi| |\xi_1|^\al$ and 
  $\jbr{\xi} \sim |\xi|$ we have
  \[
    h(\z_1,\z)
    \ls |\xi|^{1+s-\al b_2 + b'} |\xi_1|^{-2s+\frac12-\frac\al4+\al(b'+b_2)}
    \ls |\xi|^{-s+\frac32-\frac\al4+(\al+1)b'}
  \]
  where the last inequality follows from \eqref{eq:DGKP_qzbe} and
  $|\xi_1|\gs |\xi|$.
  Now $h$ is bounded because of \eqref{eq:DGKP_qebz} and $|\xi|\geq 1$.
  So let us now consider the case $|\xi|\leq 1$.
  Because of \eqref{eq:DGKP_qzbe} it follows that
  $\jbr{\la_{\max}}^{b'+b_2} \leq \jbr{\la_{\max}}^{b'+b_1}$ and therefore
  \[ 
    h(\z_1,\z) \ls
    |\xi|^{1-\lf+b'+b_1} |\xi_1|^{-2s+\frac12-\frac\al 4+\al(b'+b_1)}
  \]
  Because of \eqref{eq:DGKP_qzbz} and $|\xi|\leq 1$ we have 
  $|\xi|^{1-\lf+b'+b_1} \ls 1$.
  Because of \eqref{eq:DGKP_qzbe} we have 
  $|\xi_1|^{-2s+\frac12-\frac\al 4+\al(b'+b_1)}\ls 1$.
  So $h$ is also bounded in this case, which proves the lemma.
\end{proof}

\begin{lemma} 
  \label{le:DGKP_bilinqze}
  We have that
  \begin{equation}
    \label{eq:DGKP_bilinqze}
    \|Q_{21}(u_1,u_2)\|_{\xsb[\lf]{s-(\al+1)b_2}0{b'+b_2}} 
    \ls \|u_1\|_{\xsb{s}0b} \|u_2\|_{\xsb{s}0b}
  \end{equation}
  provided that $b>\frac12$, $b'>-\frac12$, \eqref{eq:DGKP_qebz},
  \eqref{eq:DGKP_qzbe} and \eqref{eq:DGKP_qzbz} hold.
\end{lemma}
\begin{proof}
  We have to prove that
  \begin{equation*}
    |\int_{A_{21}} k_{21}(\z_1,\z) 
      \frac{|\xi|^{-\frac12} |\xi-\xi_1|^{\frac\al 4}}
           {\jbr{\la}^b \jbr{\la_2}^b}
      f_1(\z_1)f_2(\z-\z_1)f_3(\z) d\z d\z_1|
    \ls \prod_{i=1}^3 \|f_i\|_{L^2}
  \end{equation*}
  where
  \(
    k_{21}(\z_1,\z) = \jbr{\la_1}^{-b} \jbr{\la}^{b+b'+b_2}
    \jbr{\xi}^{s-(\al+1)b_2+\lf} \jbr{\xi_1}^{-s} \jbr{\xi-\xi_1}^{-s} 
    |\xi|^{\frac32-\lf} |\xi-\xi_1|^{-\frac\al 4}.
  \)
  Now if we show that $k_{21}$ is bounded in $A_{21}$, the lemma follows
  from the refined bilinear Strichartz estimate \eqref{eq:DGKP_dbilsmdz}.
  In $A_{21}$ we have $|\la|\leq |\la_{\max}| = |\la_1|$ so that 
  \(
    \jbr{\la_1}^{-b} \jbr{\la}^{b+b'+b_2} \leq \jbr{\la_{\max}}^{b'+b_2} 
  \)
  Using that in $A_{21}$ we also have 
  $1\leq |\xi-\xi_1|\sim |\xi_1| \sim \jbr{\xi_1}$ and we get
  \begin{align*}
    k_{21}(\z_1,\z) & \ls \jbr{\la_{\max}}^{b'+b_2} 
      \jbr{\xi}^{s-(\al+1)b_2+\lf} |\xi|^{\frac32-\lf} 
      |\xi_1|^{-2s-\frac\al4} \\
      & \ls \jbr{\la_{\max}}^{b'+b_2} \jbr{\xi}^{s-(\al+1)b_2+\lf} |\xi|^{1-\lf} 
      |\xi_1|^{-2s+\frac12-\frac\al4}
  \end{align*}
  where the last inequality follows because of $|\xi|\ls |\xi_1|$
  in $A_{21}$.
  But this last expression is $h$ defined in \eqref{eq:DGKP_exprxiz}
  which was shown to be bounded in all of $\Xi_2$ in 
  Lemma~\ref{le:DGKP_bilinqzn}.
  This proves the lemma.
\end{proof}

\begin{lemma} 
  \label{le:DGKP_bilinqzz}
  We have that
  \begin{equation}
    \label{eq:DGKP_bilinqzz}
    \|Q_{22}(u_1,u_2)\|_{\xsb[\lf]{s-(\al+1)b_2}0{b'+b_2}} 
    \ls \|u_1\|_{\xsb{s}0b} \|u_2\|_{\xsb{s}0b}
  \end{equation}
  provided that $b>\frac12$, $b'>-\frac12$, \eqref{eq:DGKP_qebz},
  \eqref{eq:DGKP_qzbe} and \eqref{eq:DGKP_qzbz} hold.
\end{lemma}
\begin{proof}
  By definition \eqref{eq:DGKP_xsbndef} and duality we have to show that
  \begin{equation*}
    |\int_{A_{22}} k_{22}(\z_1,\z) 
      \frac{|\xi|^{-\frac12} |\xi_1|^{\frac\al 4}}
           {\jbr{\la}^b \jbr{\la_1}^b}
      f_1(\z_1)f_2(\z-\z_1)f_3(\z) d\z d\z_1|
    \ls \prod_{i=1}^3 \|f_i\|_{L^2}
  \end{equation*}
  where
  \(
    k_{22}(\z_1,\z) = \jbr{\la_2}^{-b} \jbr{\la}^{b+b'+b_2}
    \jbr{\xi}^{s-(\al+1)b_2+\lf} \jbr{\xi_1}^{-s} \jbr{\xi-\xi_1}^{-s} 
    |\xi|^{\frac32-\lf} |\xi_1|^{-\frac\al 4}.
  \)
  Now if we show that $k_{22}$ is bounded on $A_{22}$, the lemma follows
  from \eqref{eq:DGKP_dbilsmdd}.
  But the boundedness of $k_{22}$ follows in exactly the same way as the
  boundedness of $k_{21}$ in Lemma~\ref{le:DGKP_bilinqze}.
\end{proof}

We are now in a position to prove Theorem \ref{th:DGKP_bilest}.

\begin{proof}[Proof of Theorem \ref{th:DGKP_bilest}]
  With our definitions of $b$, $b'$, $b_1$ and $\lf$ we have 
  that $b>\frac12$, $b-b'<1$,
  \eqref{eq:DGKP_qebe}, \eqref{eq:DGKP_qebz}, \eqref{eq:DGKP_qebd}, 
  \eqref{eq:DGKP_qebv}, \eqref{eq:DGKP_qzbe} and \eqref{eq:DGKP_qzbz} hold.
  We noticed before that because of the symmetry of \eqref{eq:DGKP_nonlinest}
  in $u_1$ and $u_2$, it suffices to show 
  \eqref{eq:DGKP_nonlinest}
  for $\partial_x P_1(u_1,u_2)$ instead of $\partial_x (u_1 u_2)$ where
  $P_1$ is the operator defined in \eqref{eq:DGKP_defparapr}.
  We now decompose $\partial_x P_1(u_1,u_2)$ further as in 
  \eqref{eq:DGKP_nonldecmp}.
  Therefore we have to show for every $Q_{ij}$ that
  \begin{equation}
    \label{eq:DGKP_ineqtp}
    \|Q_{ij}(u_1,u_2)\|_{\tilde{X}} \leq C \|u_1\|_{X} \|u_2\|_{X}.
  \end{equation}
  Let us notice that by the definition of the space $\tilde{X}$ we 
  have that $\|u\|_{\tilde{X}} \leq \|u\|_{\xsb{s}00}$ and 
  $\|u\|_{\tilde{X}} 
  \leq \|u\|_{\xsb[\lf]s0{b'}} + \|u\|_{\xsb[\lf]{s-(\al+1)b_1}0{b'+b_1}}$,
  so that it suffices to control $Q_{ij}$ in one of these norms.
  The norm in $X$ is given by
  \begin{equation*}
    \|u\|_X = \inf \{ \|v\|_{X_1} + \|w\|_{X_2} \mid
    u=v+w, v\in X_1, w\in X_2 \}.
  \end{equation*}
  But by the definition of the Bourgain spaces \eqref{eq:DGKP_xsbndef} 
  and because of $b'\leq 0$ and $\lf\geq 0$ we have 
  $\|u\|_{\xsb{s}0b} \leq \|u\|_{X_1}$ and $\|u\|_{\xsb{s}0b} \leq \|u\|_{X_2}$
  which means that we have the continuous embedding $X\embeds \xsb{s}0b$.
  Therefore \eqref{eq:DGKP_ineqtp} follows from 
  \begin{equation*}
    \|Q_{ij}(u_1,u_2)\|_{\tilde{X}} \leq C \|u_1\|_{\xsb{s}0b} \|u_2\|_{\xsb{s}0b}
  \end{equation*}
  which actually holds for all of the $Q_{ij}$ except $Q_{11}$.
  Let us prove this first.
  From Lemma~\ref{le:DGKP_bilinqnn} it follows that
  \begin{align*}
    \|Q_{00}(u_1,u_2)\|_{\tilde{X}} & \leq \|Q_{00}(u_1,u_2)\|_{\xsb{s}00} \\
    & \ls \|u_1\|_{\xsb{s}0b} \|u_2\|_{\xsb{s}0b}
    \ls \|u_1\|_X \|u_2\|_X.
  \end{align*}
  In the same way, from Lemmas~\ref{le:DGKP_bilinqen},
  \ref{le:DGKP_bilinqez}, \ref{le:DGKP_bilinqzn}, \ref{le:DGKP_bilinqze}
  and \ref{le:DGKP_bilinqzz} it follows 
  \begin{align*}
    \|Q_{ij}(u_1,u_2)\|_{\tilde{X}} & \leq \|Q_{ij}(u_1,u_2)\|_{\xsb[\lf]s0{b'}} 
      + \|Q_{ij}(u_1,u_2)\|_{\xsb[\lf]{s-(\al+1)b_1}0{b'+b_1}} \\
    & \ls \|u_1\|_{\xsb{s}0b} \|u_2\|_{\xsb{s}0b}
    \ls \|u_1\|_X \|u_2\|_X
  \end{align*}
  for all of the remaining $Q_{ij}$ except $Q_{11}$.
  So it remains to consider $Q_{11}$.
  Now let us decompose $u_1\in X$ as $u_1=v_1+w_1$ with
  $v_1\in X_1$ and $w_1\in X_2$.
  For $v_1$ we have because of \eqref{eq:DGKP_bilinqeea} of 
  Lemma~\ref{le:DGKP_bilinqee}
  \begin{align*}
    \|Q_{11}(v_1,u_2)\|_{\tilde{X}} & \leq \|Q_{11}(u_1,u_2)\|_{\xsb{s}00} \\
    & \ls \|v_1\|_{\xsb{s}0{b-b'}} \|u_2\|_{\xsb{s}0b}
    \ls \|v_1\|_{X_1} \|u_2\|_X.
  \end{align*}
  For $w_1$ we have because of \eqref{eq:DGKP_bilinqeeb} of 
  Lemma~\ref{le:DGKP_bilinqee}
  \begin{align*}
    \|Q_{11}(w_1,u_2)\|_{\tilde{X}} & \leq \|Q_{11}(w_1,u_2)\|_{\xsb[\lf]s0{b'}} 
      + \|Q_{11}(w_1,u_2)\|_{\xsb[\lf]{s-(\al+1)b_1}0{b'+b_1}} \\
    & \ls \|w_1\|_{\xsb[\lf]{s-(\al+1)b_1}0{b+b_1}} \|u_2\|_{\xsb{s}0b}
    \ls \|w_1\|_{X_2} \|u_2\|_X.
  \end{align*}
  So putting these two estimates together we have
  \begin{equation*}
    \|Q_{11}(u_1,u_2)\|_{\tilde{X}} 
    \ls (\|v_1\|_{X_1} + \|w_1\|_{X_2}) \|u_2\|_X.
  \end{equation*}
  Now taking the infimum over every decomposition of $u_1$ of the
  form $u_1=v_1+w_1$ with $v_1\in X_1$ and $w_1\in X_2$ we finally
  get \eqref{eq:DGKP_ineqtp} for $Q_{11}$, which finishes the proof.
\end{proof}

For the proof of Theorem~\ref{th:DGKP_mainglob} we need the
following refined version of \eqref{eq:DGKP_nonlinest}
\
\begin{cor}
  \label{th:DGKP_rebilest}
  Let $\frac43<\al\leq 6$ and $s>\max(1-\frac34 \al,\frac14-\frac38 \al)$.
  Let $X$ and $\tilde{X}$ be defined as in Theorem~\ref{th:DGKP_bilest} 
  We then have for every $\rho>0$
  \begin{equation}
    \label{eq:DGKP_renonlinest}
    \|J_x^\rho \partial_x (u_1 u_2)\|_{\tilde{X}} \ls 
    (\|J_x^\rho u_1\|_X \|u_2\|_X + \|u_1\|_X \|J_x^\rho u_2\|_X).
  \end{equation}  
\end{cor}
\begin{proof}
  Writing
  $\partial_x (u_1 u_2)=\partial_xP_1(u_1,u_2)+\partial_xP_1(u_2,u_1)$
  where $P_1$ is the operator defined in \eqref{eq:DGKP_defparapr}
  it suffices to show
  \begin{equation}
    \label{eq:DGKP_renonlinpp}
    \|J_x^\rho \partial_x P_1(u_1,u_2)\|_{\tilde{X}} \ls
    \|u_1\|_X \|J_x^\rho u_2\|_X.
  \end{equation}
  But this follows exactly as \eqref{eq:DGKP_nonlinest},
  as the operators $J_x^\rho$ in \eqref{eq:DGKP_renonlinpp}
  only give an additional bounded term 
  $(\frac{\jbr{\xi}}{\jbr{\xi-\xi_1}})^\rho$ 
  in the dual formulations of the estimates proved in 
  Lemmas~\ref{le:DGKP_bilinqnn}-\ref{le:DGKP_bilinqzz}.
\end{proof}

\section{Proof of theorems~\ref{th:DGKP_maindgkp} and
\ref{th:DGKP_mainglob}}
\  
As the methods of proof used here are all well known, we only give
\ 
\begin{proof}[Sketch of proof of Theorem~\ref{th:DGKP_maindgkp}]
  As explained at the beginning of Section~\ref{se:DGKP_mbe}
  it suffices to consider the case $s_2=0$. Let $s=s_1$.
  For $T\leq 1$ and $u_1, u_2\in \ntf$ we define the 
  bilinear operator $\Gamma_T$ by \eqref{eq:DGKP_defgam}.
  Let $X$ and $\tilde{X}$ be defined as in Theorem~\ref{th:DGKP_bilest}
  and set $\delta:=1-(b-b')>0$.
  Then by \eqref{eq:DGKP_linestiX} and Theorem~\ref{th:DGKP_bilest}
  we have
  \begin{equation*}
    \|\Gamma_T(u_1,u_2)\|_X \ls T^\delta \|\partial_x(u_1 u_2)\|_{\tilde{X}}
    \ls T^\delta \|u_1\|_X \|u_2\|_X.
  \end{equation*}
  Therefore we can extend $\Gamma_T$ to a continuous, bilinear operator 
  $\Gamma_T : X\times X\to X$.
  As $\Gamma_T(u_1,u_2)_{|[-T,T]}$ only depends on ${u_i}_{|[-T,T]}$
  $(i=1,2)$, $\Gamma_T$ also defines a continuous, bilinear operator
  $\Gamma_T : X_T\times X_T\to X_T$.
  Furthermore by~\eqref{eq:DGKP_linestivX} we have
  \begin{equation*}
    \|\psi U_\al(\cdot) u_0\|_{X_T} \leq
    \|\psi U_\al(\cdot) u_0\|_X \ls \|u_0\|_{H^{s,0}(\R^2)}
  \end{equation*}
  for $u_0\in H^{s,0}(\R^2)$.
  So if we define
  \begin{equation}
    \label{eq:DGKP_defphi}
    \Phi_T(u,u_0) := \psi U_\al(\cdot) u_0 - \Gamma_T(u,u),
    \quad u\in X_T,\quad u_0\in H^{s,0}(\R^2)
  \end{equation}
  we have for 
  $u_0\in B_R := \{u_0\in H^{s,0}(\R^2) \mid \|u_0\|_{H^{s,0}} < R \}$
  and 
  $u, v\in \bar{A}_r := \{u\in X_T \mid \|u\|_{X_T} \leq r \}$
  that
  \begin{equation}
    \label{eq:DGKP_phie}
    \|\Phi_T(u,u_0)\|_{X_T} 
    \leq C(\|u_0\|_{H^{s,0}(\R^2)} + T^\delta \|u\|_{X_T}^2)
    \leq CR+CT^\delta r^2
  \end{equation}
  with some constant $C$ which does not depend on $R$, $r$ and $T$ and
  \begin{equation}
    \label{eq:DGKP_phiz}
    \begin{split}
      \|\Phi_T(u,u_0)-\Phi_T(v,u_0)\|_{X_T} = \|\Gamma_T(u-v,u+v)\|_{X_T} \\
      \leq C T^\delta (\|u\|_{X_T}+\|v\|_{X_T})\|u-v\|_{X_T}
      \leq 2 C T^\delta r\|u-v\|_{X_T}
    \end{split}
  \end{equation}
  So given $R>0$ we choose $r=2CR$ and 
  $T=\min\{1,(8C^2R)^{-\frac1\delta}\}$.
  Then, for fixed $u_0\in B_R$, 
  by \eqref{eq:DGKP_phie} $\Phi_T(\cdot,u_0)$ maps $\bar{A}_r$ 
  into $\bar{A}_r$
  and by \eqref{eq:DGKP_phiz} $\Phi_T(u_0,\cdot)$ is a contraction.
  By the Banach fixed point theorem there is exactly one fixed
  point of $\Phi_T$ in $\bar{A}_r$.
  Now by a well known argument
  the uniqueness of the solution $u$ follows also in $X_T$.
  Furthermore it is easy to see that the mapping
  \begin{equation*}
    \Lambda_T : X_T \times B_R \to X_T,
      \quad \Lambda_T(u,u_0) := \Phi_T(u,u_0) - u
  \end{equation*}
  is analytic. 
  Therefore a standard use of the implicit function theorem
  yields the analyticity of the flow map $F_R:u_0\mapsto u$.
\end{proof}

\begin{proof}[Sketch of proof of Theorem~\ref{th:DGKP_mainglob}]
  Let $u_0\in H^{s,0}(\R^2;\R)$ be real valued and let
  $X$ be defined as in Theorem~\ref{th:DGKP_bilest}.
  Let $T_0$ be the supremum of all $T\in (0,1]$ such that there exists 
  a unique $u\in X_T$ with $\Phi_T(u,u_0)=u$.
  We will prove that $T_0=1$.
  By Theorem~\ref{th:DGKP_maindgkp} we see that $T_0>0$.
  Let $T\in (0,T_0)$.
  Let $X^0$ be defined as $X$ in Theorem~\ref{th:DGKP_bilest},
  but with $s=0$.
  We obviously have $\|u\|_{X}=\|J^s u\|_{X^0}$.
  By \eqref{eq:DGKP_linestivX}, \eqref{eq:DGKP_linestiX} and 
  Corollary~\ref{th:DGKP_rebilest} we have
  \begin{equation*}
    \|J^s u\|_{X^0_T} \leq CR 
      + 2 C T^\delta \|u\|_{X^0_T} \|J^s u\|_{X^0_T}. 
  \end{equation*}
  Now if $T_0\leq \min(1,(8 C^2 \|u_0\|_{L^2})^{-\frac1\delta})$ we have
  by \eqref{eq:DGKP_phie} that
  \(
    \|u\|_{X^0_T} \leq 2C\|u_0\|_{L^2}
  \)
  and therefore
  \begin{equation*}
    \|J^s u\|_{X^0_T} \leq CR 
      + \frac12 \|J^s u\|_{X^0_T}.   
  \end{equation*}
  It follows that $\sup_{|t|\leq T} \|u(t)\|_{H^{s,0}} 
  \leq \tilde{C} \|J^s u\|_{X^0_T} \leq 2\tilde{C}CR$.
  As this upper bound does not depend on $T$ and
  applying Theorem~\ref{th:DGKP_maindgkp} with $u(T)$ and
  $u(-T)$ as initial values, we see that we can extend the
  solution beyond the interval $[-T_0,T_0]$.
  This contradicts the choice of $T_0$.
  Therefore we have $T_0\geq \min(1,(8 C^2 \|u_0\|_{L^2})^{-\frac1\delta})$.
  This implies that the length of the maximal interval of 
  existence does only depend on $\|u_0\|_{L^2}$.
  But the $L^2$-norm of real valued solutions $u$ of 
  \eqref{eq:DGKP_opeq} is conserved, i. e. 
  $\|u(\pm T)\|_{L^2}=\|u_0\|_{L^2}$, so if we had $T_0 < 1$,
  we could extend the solution beyond the interval $[-T_0,T_0]$ 
  which contradicts the choice of $T_0$.
\end{proof}

\bibliography{literatur}

\begin{thebibliography}{10}

\bibitem{bourgain:93_2}
J.~Bourgain.
\newblock Fourier transform restriction phenomena for certain lattice subsets
  and applications to nonlinear evolution equations. {I}. {S}chr\"odinger
  equations.
\newblock {\em Geom. Funct. Anal.}, 3(2):107--156, 1993.

\bibitem{bourgain:93_3}
J.~Bourgain.
\newblock Fourier transform restriction phenomena for certain lattice subsets
  and applications to nonlinear evolution equations. {II}. {T}he
  {K}d{V}-equation.
\newblock {\em Geom. Funct. Anal.}, 3(3):209--262, 1993.

\bibitem{bourgain:93}
J.~Bourgain.
\newblock On the {C}auchy problem for the {K}adomtsev-{P}etviashvili equation.
\newblock {\em Geom. Funct. Anal.}, 3(4):315--341, 1993.

\bibitem{ginibre:96}
Jean Ginibre.
\newblock Le probl\`eme de {C}auchy pour des {EDP} semi-lin\'eaires
  p\'eriodiques en variables d'espace (d'apr\`es {B}ourgain).
\newblock {\em Ast\'erisque}, (237):Exp.\ No.\ 796, 4, 163--187, 1996.
\newblock S\'eminaire Bourbaki, Vol.\ 1994/95.

\bibitem{iorio:98}
R.~J. I{\'o}rio, Jr. and W.~V.~L. Nunes.
\newblock On equations of {KP}-type.
\newblock {\em Proc. Roy. Soc. Edinburgh Sect. A}, 128(4):725--743, 1998.

\bibitem{isaza:06}
P.~Isaza, J.~L{\'o}pez, and J.~Mej{\'{\i}}a.
\newblock Cauchy problem for the fifth order {K}adomtsev-{P}etviashvili
  ({KPII}) equation.
\newblock {\em Commun. Pure Appl. Anal.}, 5(4):887--905, 2006.

\bibitem{isaza:01}
P.~Isaza and J.~Mej{\'{\i}}a.
\newblock Local and global {C}auchy problems for the {K}adomtsev-{P}etviashvili
  ({KP}-{II}) equation in {S}obolev spaces of negative indices.
\newblock {\em Comm. Partial Differential Equations}, 26(5-6):1027--1054, 2001.

\bibitem{kadomtsev:70}
B.B. Kadomtsev and V.I. Petviashvili.
\newblock {On the stability of solitary waves in weakly dispersing media}.
\newblock {\em Sov. Phys., Dokl.}, 15:539--541, 1970.

\bibitem{kenig:91}
C.~E. Kenig, G.~Ponce, and L.~Vega.
\newblock Oscillatory integrals and regularity of dispersive equations.
\newblock {\em Indiana Univ. Math. J.}, 40(1):33--69, 1991.

\bibitem{saut:06}
L.~Molinet, J.~C. Saut, and N.~Tzvetkov.
\newblock Remarks on the mass constraint for {KP} type equations.
\newblock 2006, arXiv:math.AP/0603303.

\bibitem{saut:93}
J.~C. Saut.
\newblock Remarks on the generalized {K}adomtsev-{P}etviashvili equations.
\newblock {\em Indiana Univ. Math. J.}, 42(3):1011--1026, 1993.

\bibitem{saut:99}
J.~C. Saut and N.~Tzvetkov.
\newblock The {C}auchy problem for higher-order {KP} equations.
\newblock {\em J. Differential Equations}, 153(1):196--222, 1999.

\bibitem{saut:00}
J.~C. Saut and N.~Tzvetkov.
\newblock The {C}auchy problem for the fifth order {KP} equations.
\newblock {\em J. Math. Pures Appl. (9)}, 79(4):307--338, 2000.

\bibitem{takaoka:00}
H.~Takaoka.
\newblock Global well-posedness for the {K}adomtsev-{P}etviashvili {II}
  equation.
\newblock {\em Discrete Contin. Dynam. Systems}, 6(2):483--499, 2000.

\bibitem{takaoka:01_2}
H.~Takaoka.
\newblock Well-posedness for the {K}adomtsev-{P}etviashvili {II} equation.
\newblock {\em Adv. Differential Equations}, 5(10-12):1421--1443, 2000.

\bibitem{takaoka:01}
H.~Takaoka and N.~Tzvetkov.
\newblock On the local regularity of the {K}adomtsev-{P}etviashvili-{II}
  equation.
\newblock {\em Internat. Math. Res. Notices}, (2):77--114, 2001.

\bibitem{tzvetkov:99}
N.~Tzvetkov.
\newblock On the {C}auchy problem for {K}adomtsev-{P}etviashvili equation.
\newblock {\em Comm. Partial Differential Equations}, 24(7-8):1367--1397, 1999.

\bibitem{tzvetkov:00}
N.~Tzvetkov.
\newblock Global low-regularity solutions for {K}adomtsev-{P}etviashvili
  equation.
\newblock {\em Differential Integral Equations}, 13(10-12):1289--1320, 2000.

\end{thebibliography}

\end{document}